\documentclass[12pt]{article}
\usepackage{amssymb}
\usepackage{epsfig}

\newcommand{\triangleboard}{\textsf{Triangle}}

\newenvironment{packed_enumerate}{
\setlength{\parsep}{0pt}
\setlength{\parskip}{0pt}
\begin{enumerate}
  \setlength{\itemsep}{1pt}
  \setlength{\parsep}{0pt}
  \setlength{\parskip}{0pt}
}{\end{enumerate}}

\begin{document} 

\begin{center}
\vspace*{0.1in}
{\huge Triangular Peg Solitaire Unlimited}
\vskip 20pt
{\bf George I. Bell}\\
December 2004\footnote{Original version at {\tt http://gpj.connectfree.co.uk/gpjr.htm} \\
\hspace*{0.25in}Converted to \LaTeX ~by the author with some modifications to the text, November 2007.} \\
{\tt gibell@comcast.net}\\
\end{center}
\vskip 30pt 
\centerline{\bf Abstract}
\noindent
Triangular peg solitaire is a well-known one-person game or puzzle.
When one peg captures many pegs consecutively,
this is called a sweep.
We investigate whether the game can end in a dramatic fashion,
with one peg sweeping all remaining pegs off the board.
For triangular boards of side 6 and 8
(with 21 and 36 holes, respectively)
the geometrically longest sweep can occur as the final move in a game.
On larger triangular boards, we demonstrate how to construct solutions
that finish with arbitrarily long sweeps.
We also consider the problem of finding solutions
that minimize the total number of moves
(where a move is one or more consecutive jumps by the same peg).

\pagestyle{myheadings}
\markright{The Games and Puzzles Journal---Issue 36, November-December 2004\hfill}

\thispagestyle{empty} 
\baselineskip=15pt 
\vskip 30pt 

\section{Introduction}

\noindent
Peg solitaire on a 15-hole triangular board is an
old puzzle but it remains popular.
In the United States, one can find these puzzles
at tables in $\mbox{Cracker Barrel}\textsuperscript{\textregistered}$~restaurants.
The 15-hole puzzle is also amenable to exhaustive computer search
and this is a common programming assignment for computer science classes.

\begin{figure}[htb]
\centering
\epsfig{file=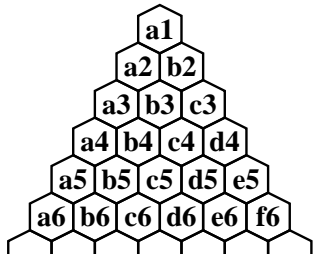}
\caption{The general triangular board with hole coordinates.}
\label{fig1}
\end{figure}

\noindent
Here we consider peg solitaire on an (equilateral) triangular board with
$n$ holes on each side.
This board will be referred to as \triangleboard$(n)$
and can be conveniently presented on an array
of hexagons (Figure~\ref{fig1}).
The \triangleboard$(n)$ board has $T(n)=n(n+1)/2$ holes,
where $T(n)$ is the $n$'th triangular number.
The notation used to identify the holes in the board is shown in
Figure~\ref{fig1}---it differs from that given in Beasley \cite{Beasley}
and is based on the system for square lattice boards.
A nice property of this notation is that the top corner hole is always ``a1".

\noindent
The rules of the game are simple: start from a
board with a peg in every hole but one,
called the \textbf{starting vacancy}.
Then jump one peg over another into an empty hole,
removing the jumped peg from the board.
The goal is to choose a sequence of jumps to finish with one peg
(the coordinate of this final peg is called
the \textbf{finishing location}).
This general problem of going from a board position with one peg missing to
a board position with one peg will be referred to as a
\textbf{peg solitaire problem}.
The special case where the starting vacancy and finishing location
are the same is referred to as a \textbf{complement problem}.
GPJ \#28 \cite{GPJ28} considered peg solitaire on boards based on a square lattice.
Triangular peg solitaire differs in that the holes are based on a triangular
lattice, with a maximum of six jumps available into any hole
(rather than four).
However the theory of the game (in particular the fundamental classes)
is very similar, see \cite{Beasley}, or \cite{Bellweb} for details.

\noindent
Any solution to a problem on the 15-hole triangular board consists
of exactly 13 \textbf{jumps} because we start with 14 pegs and
finish with 1.
However, when the same peg jumps one or more pegs consecutively,
we call this one \textbf{move}.
Given a particular peg solitaire problem,
what is the solution with the least number of moves?
While this question has historically played an important role in
peg solitaire on the standard 33-hole cross-shaped board,
it has barely been considered for triangular boards.
The following terminology is used in referring to moves involving multiple
jumps: when a peg removes $i$ pegs in a single move, we refer to it
as a \textbf{sweep}, or more specifically, an \textbf{\textit{i}-sweep}.

\noindent
After attempting a peg solitaire problem,
many people get the idea to try to solve it
\textit{backwards} from the final peg.
What is not so obvious is that this is \textit{exactly the same}
as the original game, where the concepts of ``hole" and ``peg" are interchanged.
In fact the solution to any peg solitaire problem
really contains \textit{two} solutions:
the original (``forward" solution) plus this ``backward" solution,
where the individual jumps are executed in the same direction, but in reverse
order,
and the starting vacancy and finishing location are swapped.
An important observation is that when an $i$-sweep is reversed,
we must reverse the individual jumps, and it becomes $i$ single
jump moves.
In other words, sweeps in forward solutions \textit{cannot be} sweeps in
reversed solutions (and vice versa).

\noindent
Backward play is hard to comprehend, because our brain does not easily
interchange the concepts of ``hole" and ``peg".
It is easier to understand ``backward play" by realizing that it is the
same as forward play from the \textbf{complement} of the current board position
(the complement of a board position is where we replace every hole by a peg,
and vice versa).
This leads to a theorem used extensively in the remainder of this paper.
Suppose we have a board position $\mathcal B$ and we wonder if it could
appear during a solution to a peg solitaire problem.

\noindent
\textbf{Forward/Backward Theorem}: Board position $\mathcal B$ can
appear during a solution to a peg solitaire problem if and only if both
\begin{packed_enumerate}
\item $\mathcal B$ can be reduced to a single peg
(the finishing location) using peg solitaire jumps.
\item The complement of $\mathcal B$ can be reduced to a single peg
(the starting location) using peg solitaire jumps.
\end{packed_enumerate}
A proof of this theorem just boils down to the equivalence of
playing the game forward or backward.
It may seem obvious, but is key to understanding
how to reach complex sweep positions.

\noindent
A practical problem that arises is finding a triangular board to play on.
Fifteen hole \triangleboard$(5)$ boards are common, but cannot be
easily extended
to larger triangular boards.
The best board I have found is a Chinese Checkers set, which allows one to
play on boards as large as \triangleboard$(13)$.
To help follow the arguments below, I recommend
playing out the solutions on a Chinese Checkers set.
It is particularly helpful to see solutions played forward and backward.

\section{Maximal sweeps on odd-side triangular boards}

\noindent
Triangular boards support the longest sweeps
of any peg solitaire board.
This is because from any board location the total number
of possible jumps \textit{is even}.
If the board size $n$ is odd, there exist sweeps that jump into
or over every single location on the board.
Such sweeps are \textbf{maximal} in the sense that they are
the longest sweep geometrically possible on the board.
The figure below shows examples of the first four maximal sweeps

\begin{figure}[htb]
\centering
\epsfig{file=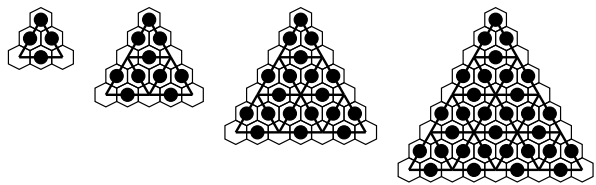}
\caption{Maximal Sweeps on \triangleboard$(n)$, where $n=3$, $5$, $7$ and $9$.
These sweeps are shown starting and ending at a1,
but can begin and end at other board locations.}
\label{fig2}
\end{figure}

\noindent
The length of this sweep is $3T((n-1)/2)=3(n^2-1)/8$.
Related geometrical tours on a triangular lattice have appeared in
GPJ \#20 \cite{GPJ20} under ``Trapezoidal Tours".
As these sweep patterns become larger and larger, how can we
be sure a peg can traverse the entire pattern?
If we consider the pattern of the sweep as a graph,
one of the most elementary theorems in graph theory, due to Euler,
ensures that such a traversal is always possible when there are at most two
nodes of odd degree.  Since maximal sweeps have all nodes of even degree,
they can always be traversed,
and the starting and ending board locations must be the same.
This theorem also guarantees that the more complex sweep patterns seen later
can be traversed.

\noindent
Note that the sweep length divided by the board size
approaches $3/4$ as the board size increases.
If this sweep is the final move to a solitaire game it removes
nearly $3/4$ of the pegs that we started with in one move!
But can these maximal sweeps be realized during a solitaire game?
If we take the complement of the sweep pattern, not a single jump is possible.
Thus, by the forward-backward theorem, maximal sweep patterns cannot appear in
solitaire games on odd size boards.

\section{Maximal sweeps on even-side triangular boards}

\noindent
The same maximal sweep on an odd triangular board is
also maximal on the even board one size larger,
because the added row cannot be reached to extend the sweep.
However, the added row can make it possible to reach the sweep position during
a peg solitaire game.
On the \triangleboard$(6)$ and \triangleboard$(8)$ boards this can
be worked out by hand,
Figure~\ref{fig3} shows this process on \triangleboard$(6)$:

\begin{figure}[htb]
\centering
\epsfig{file=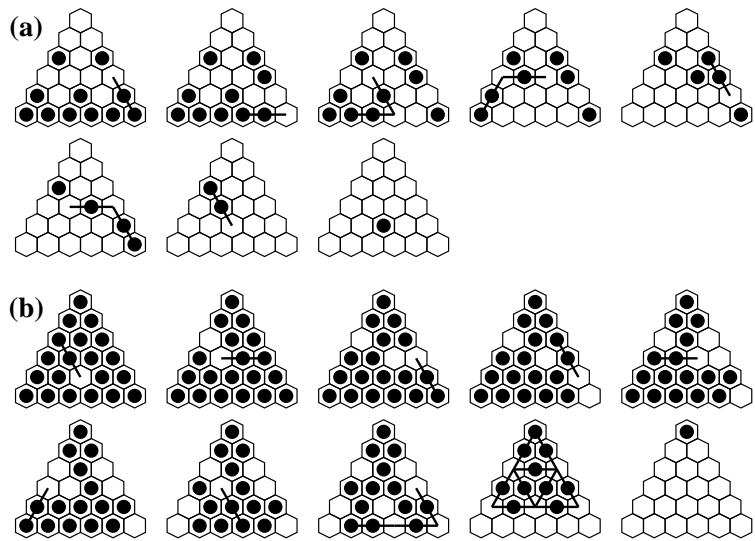}
\caption{Constructing a solitaire solution that finishes with the 9-sweep to $a1$.
(a) Backward: playing from the complement of the sweep pattern to c5.
(b) Forward: The completed solution from the c5 vacancy ending with a 9-sweep to a1.}
\label{fig3}
\end{figure}

\noindent
The forward solution in Figure~\ref{fig3}b has only 9 moves.
In the final section of this paper, we
prove that it's \textit{impossible} to solve
this problem in fewer than 9 moves.
Thus, besides containing a maximal length sweep,
the solution of Figure~\ref{fig3}b solves the problem
in the minimum number of moves.
This solution was discovered before 1975 by
Harry O. Davis \cite{Gardner}.

\noindent
There are three problems on \triangleboard$(6)$ that can contain a 9-sweep:
\begin{packed_enumerate}
\item Vacate c5, play to finish at a1 with the last move a 9-sweep
(Figure~\ref{fig3}).
\item Vacate c5, play to finish at a4 with the last move a 9-sweep.
\item Vacate c5, play to finish at a4 with the second to the
last move a 9-sweep.
\end{packed_enumerate}
The reader is invited to solve problems 2 and 3.
As in Figure~\ref{fig3},
the trick is to figure out what the final sweep must look like,
and then solve backwards
(playing forward from the complement of the sweep position).
Then exactly reverse the jumps in your ``backward" solution and you will
reproduce the sweep position.

\begin{figure}[htb]
\centering
\epsfig{file=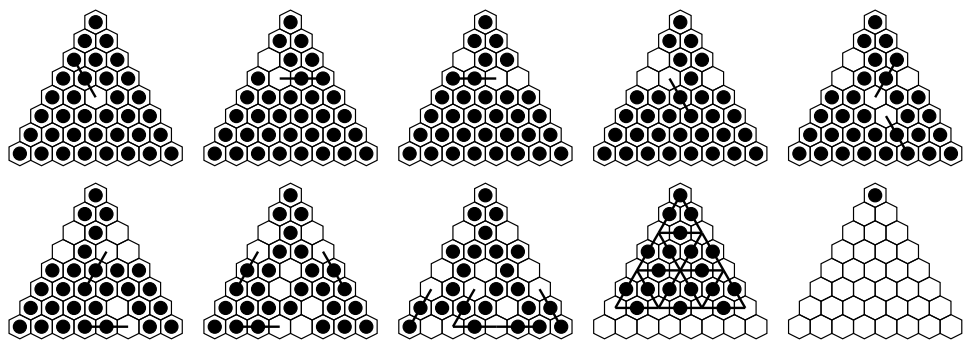}
\caption{A 15-move solution to problem \#1 on \triangleboard$(8)$ (finishing
with an 18-sweep).  Note: more than one move is sometimes shown between board snapshots.}
\label{fig4}
\end{figure}

\noindent
On \triangleboard$(8)$, the 18-sweep can also be reached.
There are three problems on this board that can contain the 18-sweep:
\begin{packed_enumerate}
\item Vacate c5, play to finish at a1 with the last move an 18-sweep
(Figure~\ref{fig4}).
\item Vacate b6 or e6, play to finish at c8 with the last move an 18-sweep.
\item Vacate c5, play to finish at b6 with the
second to the last move an 18-sweep.
\end{packed_enumerate}
These problems are more difficult than those on \triangleboard$(6)$,
but are still quite reasonable to work out by hand.
Figure~\ref{fig4} shows the solution to the first problem, played forward.
This solution was discovered, as usual, by attempting to play backward from
the complement of the sweep position.
This solution contains 15-moves, but it is possible to solve this problem
in 14~moves (without the 18-sweep), so this solution does not have
the minimum number of moves.

\noindent
On \triangleboard$(10)$, a computational search for a peg solitaire
solution containing a maximal 30-sweep has come up empty
(although a solution was found which ends with a 29-sweep).
\triangleboard$(12)$ does not appear to have a
maximal 45-sweep solution either.
I haven't checked all possible configurations of this 45-sweep,
but my program has shown the most obvious candidates can't
be reached from a single vacancy start.
It appears that the \triangleboard$(8)$
board is the largest triangular board for which a solution
to a peg solitaire problem can contain a maximal sweep.

\section{Long sweeps on arbitrarily large boards}

\begin{figure}[htb]
\centering
\epsfig{file=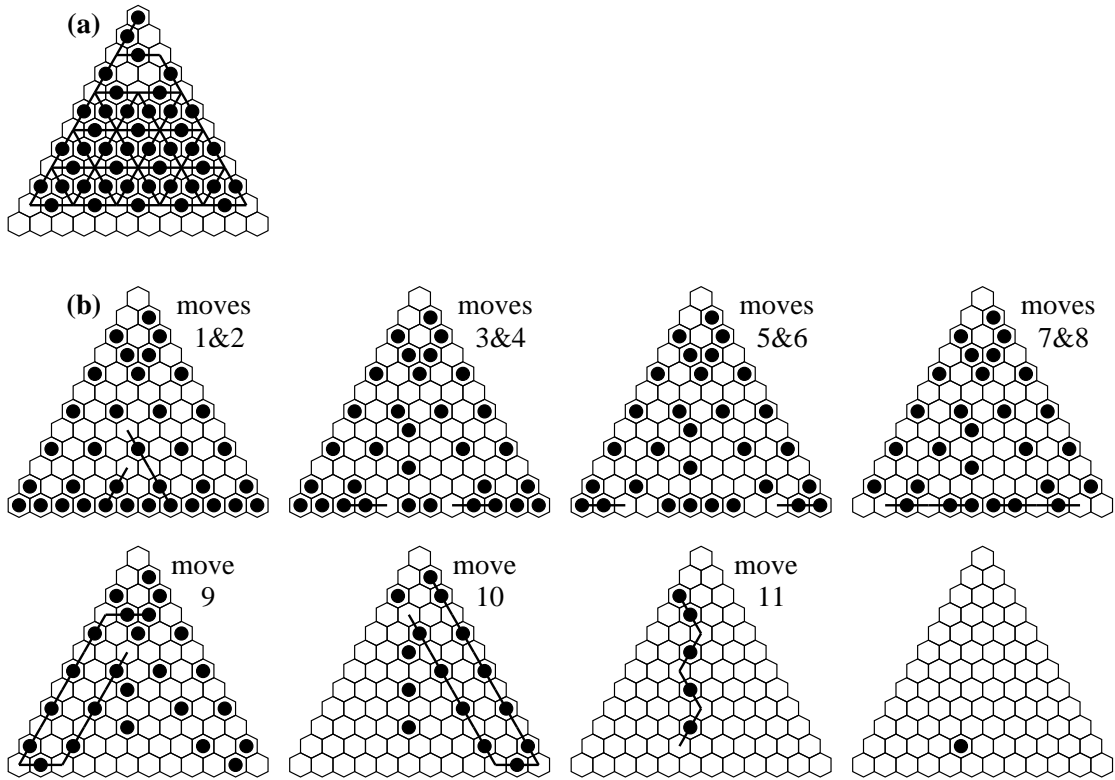}
\caption{Building a solitaire solution that finishes with the 42-sweep
(a) Forward: the finishing 42-sweep (from a1 to a3)
(b) Backward: showing how to reduce the complement of the sweep
pattern to one peg.}
\label{fig5}
\end{figure}

\noindent
Although \textit{maximal} sweeps appear not to be attainable
in peg solitaire problems on large triangular boards,
it turns out sweeps of only \textit{slightly reduced} length are possible.
A very special solution was discovered by hand on the 78-hole
\triangleboard$(12)$ board.
This solution finishes with a 42-sweep when run forward,
and is shown in Figure~\ref{fig5}.

\noindent
This solution is remarkable because it can be extended
to even larger triangular boards.
We can extend the final sweep to cover the bottom of
\triangleboard$(14)$ and the
complement can still be reduced to a single peg.
We do this by keeping moves along the bottom row the same,
while extending other moves vertically.
In particular moves 1 \& 2 end at the same board locations but begin
from the bottom row of the board.
The U-shaped moves 9 \& 10 become vertically
elongated, but have the same starting and ending board locations.
Finally additional moves need to be added after move 11 to
reduce the remaining pattern of pegs to a single survivor
(which does not end up at the same location
as in Figure~\ref{fig5}b).
We can continue stepping the board size up by 2 and still reduce the
complement to a single peg up to \triangleboard$(20)$.
However if we try this on \triangleboard$(22)$, we find we can no
longer reduce the remaining pattern to a single peg.

\noindent
As the board is extended, the middle of the lower portion of the board
becomes similar to the \textit{original problem} on \triangleboard$(12)$.
This key insight suggests that we may be able to perform an inductive step,
and \textit{extend the board indefinitely}.
This is in fact possible, and will be described below.
In effect we can construct solutions to peg solitaire problems on
triangular boards of arbitrarily large size.
Not only that, these solutions finish with arbitrarily long sweep moves!

\begin{figure}[htb]
\centering
\epsfig{file=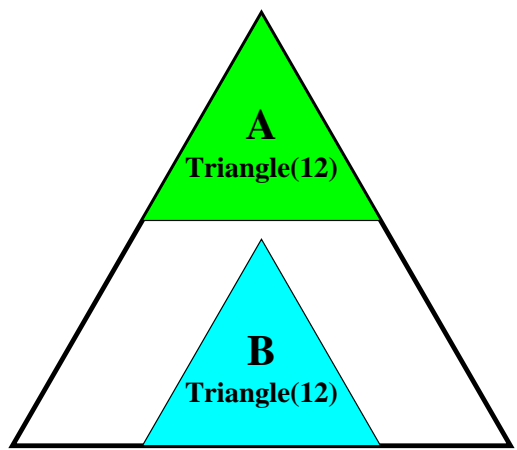}
\caption{Fitting the inductive components together to make a long
finishing sweep on \triangleboard$(24)$.}
\label{fig6}
\end{figure}
\begin{figure}[htb]
\centering
\epsfig{file=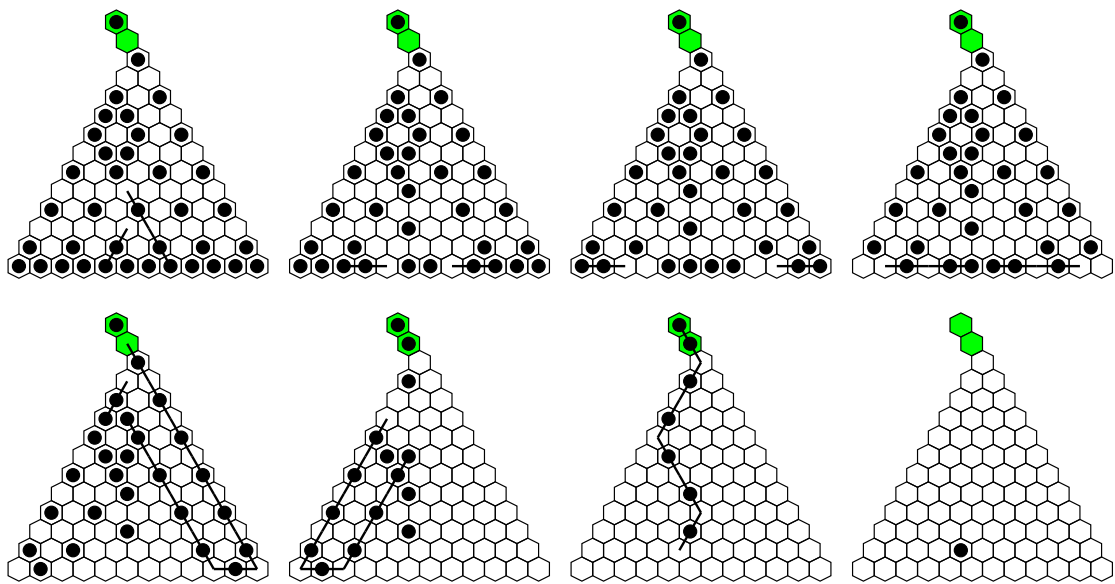}
\caption{Inductive Component $B$.  Backward: showing how to reduce the
complement of the sweep pattern to one peg.
The top two holes (in green) are part of the board above.}
\label{fig7}
\end{figure}

\noindent
To complete the inductive step, we add another solution for
\triangleboard$(12)$ underneath the first one to obtain a
solution on the 300-hole \triangleboard$(24)$ board.
Figure~\ref{fig6} shows the overall geometry of the solution; we use
``component $A$" (Figure~\ref{fig5}) in the upper part of the board and
``component $B$" (Figure~\ref{fig7}) in the lower part of the board.
In the reversed solution, the remaining (white) areas of the board
are cleared by extensions of the moves to clear $A$, and the final move
snakes down vertically through both
$A$ and $B$ to the bottom of the board.

\noindent
The final sweep pattern will be close to the maximal sweep,
but component $A$ contains a ``defect" in this sweep pattern, and so
must component $B$.  Figuring out exactly where to place this defect
in component $B$ is critical to making everything work out.
At first I tried putting the defect symmetrically over the $y$-axis
of the board, but this never seems to work out.
What did work was to put it slightly off center.

\noindent
While component $A$ can be considered a solution itself on
\triangleboard$(12)$,
component $B$ is inherently tied to the board above it, for there must be
moves which link the two boards.
Figure~\ref{fig7} shows how to solve component $B$,
which looks very similar to component $A$.
Note, however, that the U-shaped moves actually go in
the opposite direction from before.
The entire solution is constructed to enable the final move to pass
vertically down the board.

\noindent
The final sweep pattern always starts at a1 and finishes at a3.
The initial vacancy, or finishing location for the reversed solution is
always directly under a3 on the lowest row that has a hole in
vertical alignment with a3 (on \triangleboard$(24)$, this starting
vacancy is at k23).
The entire solution on \triangleboard$(24)$, when executed in the forward
direction,
has a final move which is a 191-sweep (Figure~\ref{fig8}):
\begin{figure}[htb]
\centering
\epsfig{file=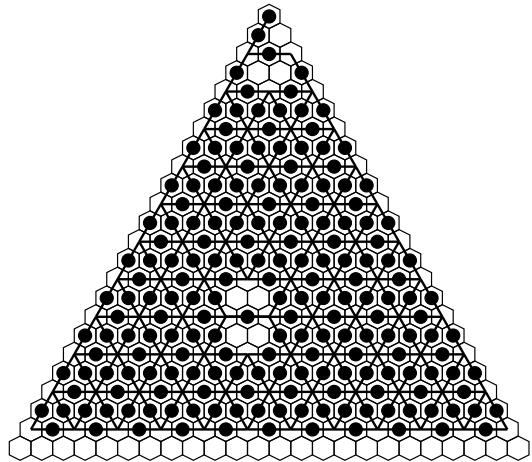}
\caption{Forward: the finishing 191-sweep on Triangle(24).
This sweep pattern is the maximal 198-sweep with two ``defects".}
\label{fig8}
\end{figure}

\noindent
By stacking additional copies of component $B$ under
the diagram of Figure~\ref{fig6},
we can extend this process indefinitely.
The net result is that on \triangleboard$(12i)$,
we can construct a solution to a solitaire problem
that finishes with a sweep of length $54i^2-13i+1$
(this is $4i-1$ shorter than the the maximal sweep length).
The board itself has $72i^2+6i$ holes, so asymptotically
this finishing sweep removes $3/4$ of the pegs on the board.
The forward solution consists of $18i^2+19i-3$
jumps, followed by the final sweep move.
By reordering a few jumps,
we can reduce the total number of moves in the
forward solution slightly to $18i^2+14i-3$.
Table~\ref{table1} gives statistics on these
solutions as the board size increases.
\begin{table}[htbp]
\begin{center} 
\begin{tabular}{ | l | c | c | c | c | c | }
\hline
 & & Board size & Final sweep & \% pegs removed & Forward solution \\		
$i$ & Board & (\# holes) & length (\# jumps) & by final sweep & length (\# moves) \\
\hline
\hline
$1$ & \triangleboard$(12)$ & 78 & 42 & 55.3\% & 29 \\
\hline
$2$ & \triangleboard$(24)$ & 300 & 191 & 64.1\% & 97 \\
\hline
$3$ & \triangleboard$(36)$ & 666 & 448 & 67.5\% & 201 \\
\hline
$10$ & \triangleboard$(120)$ & 7,260 & 5,271 & 72.6\% & 1,937 \\
\hline
\end{tabular}
\caption{Statistics on long sweep solutions on \triangleboard$(12i)$.} 
\label{table1}
\end{center} 
\end{table}

\noindent
Although this construction technique has given us solutions on boards with sides
a multiple of 12, it is not hard to extend it to \triangleboard$(n)$
for any \textit{even} $n\ge 12$.
The way to do this is using the same technique we used to
extend component $A$ on \triangleboard$(12)$ to
\triangleboard$(14)$, \triangleboard$(16)$, etc.

\section{Short solutions}

Given a board and a (solvable) peg solitaire problem,
what is the least number of moves that can solve it?
While this question has been important in the history of the standard 33-hole
cross-shaped board, it has not received much attention on triangular boards.
Informally, this question can be rephrased: ``What is the solution which involves
touching the smallest number of pegs?"

\noindent
In 1966, Harry O. Davis studied short solutions on the \triangleboard$(5)$
board analytically \cite{Davis66}.
He was able to find minimal solutions ``for all starting locations"
and prove that his solutions were the shortest possible.
In particular, he found a 10 move solution to the a1-complement,
and proved that the problem could not be solved in less than 10 moves.

\noindent
I have now completed exhaustive computer calculations on all peg solitaire
problems on boards up to \triangleboard$(7)$, and
all complement problems on \triangleboard$(8)$.
Table~\ref{table2} summarizes my results,
for more information with examples and complete lists of
shortest solutions see \cite{Bellweb}.

\begin{table}[htbp]
\begin{center} 
\begin{tabular}{ | l | c | c | c | c | c | c | c | c | c | c | }
\hline
 & Total & \multicolumn{2}{| c |}{Peg solitaire problems} &
\multicolumn{7}{| c |}{\# with minimal solution length} \\		
Board & holes & total & solvable & 9 & 10 & 11 & 12 & 13 & 14 & 15 \\
\hline
\hline
\triangleboard$(5)$ & 15 & 17 & 12 & 2 & 6 & 4 & - & - & - & - \\
\hline
\triangleboard$(6)$ & 21 & 29 & 29 & 16 & 11 & 2 & - & - & - & - \\
\hline
\triangleboard$(7)$ & 28 & 27 & 27 & - & - & - & 19 & 8 & - & - \\
\hline
\triangleboard$(8)$ & 36 & 80 & 80 & - & - & - & - & 1$\dagger$ & 5$\dagger$ & 2$\dagger$ \\
\hline
\end{tabular}
\caption{Summary of minimal solution lengths (found by computational search)
on triangular boards of side $5$--$8$
($\dagger$ - complement problems only).} 
\label{table2}
\end{center} 
\end{table}

\noindent
In Table~\ref{table2}, we count only \textit{distinct} peg solitaire
problems that cannot be reduced to another problem by means
of rotation or reflection of the board.
For example, on the \triangleboard$(5)$ board,
Table~\ref{table2} indicates that there are 17 distinct peg
solitaire problems, but only 12 are solvable.
Of these 12 solvable problems,
two can be done in a minimum of 9 moves,
six in 10 moves and four in 11 moves.
Surprisingly, over half the problems on the
\triangleboard$(6)$ board can be solved in 9 moves,
so on average, it's possible to solve problems on
this board in \textit{fewer} moves than for \triangleboard$(5)$!

\noindent
There are 80 different problems on \triangleboard$(8)$, I have only
run the 8 complement problems plus a few non-complement problems.
Among the complement problems,
only one can be accomplished in 13 moves,
the a7-complement (Figure~9)
[due to board symmetry, this problem is equivalent
to the complement problem at 5 other board locations: a2, b2, b8, g7 or g8].
The computer run also determined that this 13-move solution is \textit{unique}
in the sense that any other 13-move solution to the a7-complement must contain
the \textit{exact same set of jumps}, although not necessarily in the same order.
I have found several other 13-move solutions to \textit{non-complement} problems
on the \triangleboard$(8)$ board.

\begin{figure}[htb]
\centering
\epsfig{file=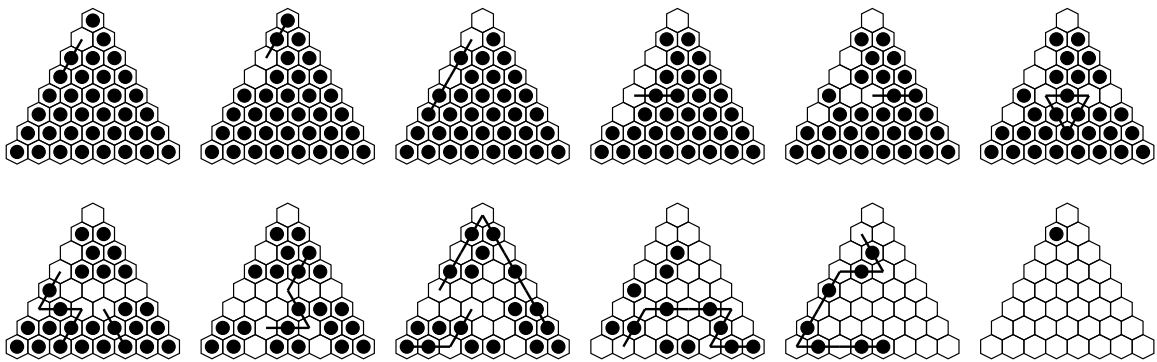}
\caption{The minimal 13-move solution to the a7-complement on \triangleboard$(8)$.
Note: more than one move is sometimes shown between board snapshots.}
\label{fig9}
\end{figure}

\section{Bounds on the length of the shortest solution}

\begin{figure}[htb]
\centering
\epsfig{file=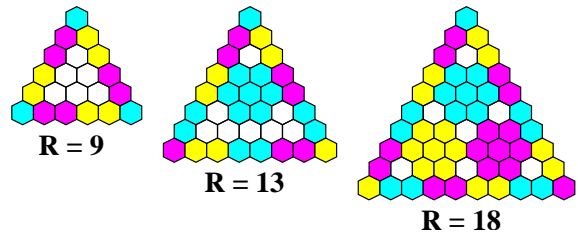}
\caption{Dividing \triangleboard$(6)$, \triangleboard$(8)$ and
\triangleboard$(10)$ into ``Merson Regions"
($R$ is the number of regions).}
\label{fig10}
\end{figure}

\noindent
It is quite difficult to determine the shortest solution
on boards larger than \triangleboard$(8)$.
However, we can obtain a lower bound on the length of the
shortest solution by dividing the board into ``Merson Regions"
(named after Robin Merson who first used this concept
in 1962 \cite[p. 203]{Beasley}).
The shape of a region is chosen such that when it is completely
filled with pegs, there is no way to remove a peg in the region
without a move that originates in the region.
On the edge of the board the regions can be corners or two consecutive holes,
but in the interior we must use hexagons 2 holes
on a side (7 holes total).
Figure~\ref{fig10} shows three triangular boards divided into Merson Regions.

\noindent
Any region that starts out full must have at least
one move starting from inside it.
Since the starting position has every hole filled
by a peg except one,
all regions start full except possibly the region that
contains the starting location.
If we start in a corner, then this region starts out
empty but is filled by the first move,
hence there still must be a move out of this corner region.
We can summarize these results as follows:
If $R$ is the number of Merson Regions, then
\begin{packed_enumerate}
\item If the starting vacancy is a corner,
or is not in any region
(white board locations in Figure~\ref{fig10}),
then any solution to the peg solitaire problem
(no matter where it finishes) has at least $R$ moves.
\item If the starting vacancy is in a region (but not a corner),
then any solution to the peg solitaire problem
(no matter where it finishes) has at least $R-1$ moves.
\end{packed_enumerate}
For example, on \triangleboard$(6)$, this proves that any solution to
``Vacate c5, finish at a1" will take at least 9 moves.
On \triangleboard$(8)$, this analysis indicates that the a7-complement
must take at least 12 moves.
It can be seen that the 13-move solution in
Figure~\ref{fig9} contains one ``extra" move above this minimum,
in that there are two moves out of the central (blue) hexagon region.
Therefore, 13 moves is not proved the shortest possible by this method,
although exhaustive computer search indicates no 12-move solution exists.

\noindent
As the triangular board becomes very large,
the number of interior hexagons eventually dominates the number of regions,
because this is the only region count growing quadratically.
We can then ``tile" the board with these hexagons without leaving any gaps
(except near the edge of the board).
So no solution can be shorter than the number of
holes in the board divided by 7.
No matter how large the triangular board,
we cannot hope to find a solution which does better
than a 7-sweep averaged over all moves.
Of course, averaging anywhere near this would be quite
remarkable---note that the 13-move solution in Figure~\ref{fig9}
removes an average of only $34/13\approx 2.6$ pegs per move.

\noindent
Using the long sweep solution in the last section,
we can obtain an upper bound on the length
of the shortest solution of any problem on \triangleboard$(12i)$.
This upper bound is a solution of length $18i^2+14i-3$
on a board of size $72i^2+6i$.
As $i$ increases, this shows that there exist solutions that
average (asymptotically) a 4-sweep for every move.
Therefore, combining the result of the last paragraph,
we see that for large triangular boards,
the shortest solution must average between
4 and 7 pegs captured per move.
Of course, finding any such solution will be very difficult.

\end{document}